\documentclass[oneside,a4paper,reqno,11pt]{amsart}
\usepackage{latexsym,amsfonts,amssymb,amsmath,amscd}
\usepackage{graphicx,color,ifpdf}
\allowdisplaybreaks

\usepackage[bookmarks=true,pdfstartview=FitH, pdfborder={0 0 0},
colorlinks=true,citecolor=blue, linkcolor=blue]{hyperref}

\usepackage{aliascnt}

\newtheorem*{thma}{Theorem A}
\newtheorem*{thmb}{Theorem B}
\newtheorem*{thmc}{Theorem C}

\newtheorem{theorem}{Theorem}[section]
\newaliascnt{lemma}{theorem}
\newtheorem{lemma}[lemma]{Lemma}
\aliascntresetthe{lemma}

\newaliascnt{proposition}{theorem}
\newtheorem{proposition}[proposition]{Proposition}
\aliascntresetthe{proposition}

\newaliascnt{corollary}{theorem}
\newtheorem{corollary}[corollary]{Corollary}
\aliascntresetthe{corollary}

\theoremstyle{definition}
\newaliascnt{definition}{theorem}

\aliascntresetthe{definition}

\newtheorem*{acknowledgements}{Acknowledgements}

\newaliascnt{example}{theorem}

\aliascntresetthe{example}

\newaliascnt{exercise}{theorem}

\aliascntresetthe{exercise}

\newaliascnt{question}{theorem}

\aliascntresetthe{question}

\newaliascnt{problem}{theorem}

\aliascntresetthe{problem}

\theoremstyle{remark}
\newaliascnt{remark}{theorem}
\newtheorem{remark}[remark]{Remark}
\aliascntresetthe{remark}

\newaliascnt{notation}{theorem}

\aliascntresetthe{notation}

\newaliascnt{fact}{theorem}

\aliascntresetthe{fact}

\numberwithin{equation}{theorem}%
\numberwithin{figure}{theorem}

\renewcommand{\Bbb}[1]{\mathbb{#1}}
\newcommand{\abs}[1]{\left|#1\right|}

\newcommand{\injrad}{\operatorname{injrad}}

\newcommand{\dist}{\operatorname{dist}}
\newcommand{\curv}{\operatorname{curv}}

\begin{document}
\title[Homotopy lifting property]{Homotopy lifting property of an
$e^\epsilon$-Lipschitz and co-Lipschitz map
}
\author{Shicheng Xu}
\email{shichengxu@gmail.com}
\address{Department of Mathematics, Nanjing University, Nanjing,
China}

\thanks{Project 11171143 supported by NSFC}

\date{\today}

\begin{abstract}
An $e^\epsilon$-Lipschitz and co-Lipschitz map, as a metric analogue
of an $\epsilon$-Riemannian submersion, naturally arises from a
sequence of Alexandrov spaces with curvature uniformly bounded below
that converges to a space of only weak singularities. In this paper we
prove its homotopy lifting property and its homotopy stability in
Gromov-Hausdorff topology. 
\end{abstract}
\maketitle

A map $f:X\to Y$ between two metric spaces is called an {\it
$e^\epsilon$-Lipschitz and co-Lipschitz} \cite{Ka07, RX12} (briefly,
$e^\epsilon$-LcL), if for any $p\in X$, and any $r>0$, the
metric balls satisfy $$B_{e^{-\epsilon}r}(f(p))\subseteq
f(B_r(p))\subseteq B_{e^\epsilon r}(f(p)).$$
A $1$-LcL preserves metric balls exactly and is called a {\it
submetry}.
An $e^\epsilon$-LcL naturally arises from a
sequence of Alexandrov spaces with curvature uniformly bounded below
that converges to a space of only weak singularities.
Recall that an Alexandrov space with curvature bounded below, curv
$\ge\kappa$, \cite{BGP92} is a complete length metric space such that
any geodesic triangle looks fatter than that with the same
side-lengths in the simply-connected space of constant curvature
$\kappa$. The Hausdorff dimension of an Alexandrov space is always an
integer or $\infty$.  

In this notes we will study the fibration
arising from of an $e^\epsilon$-LcL (Theorem A and Theorem B), and
give some applications on convergence of Alexandrov spaces in
Gromov-Hausdorff topology (see Section 2).
The first result is about the homotopy lifting
property of an $e^\epsilon$-Lipschitz and co-Lipschitz map.

\begin{thma}\label{thm:A}
A proper $(1.023)$-LcL $f:X\to B$ from a finite-dimensional Alexandrov
space $X$ with curvature bounded below to a $n$-dimensional
Riemannian manifold $B$ is a Hurewicz fibration i.e., satisfying the
homotopy lifting property with respect to any space.
\end{thma}

\begin{remark}\label{rem:thma}
In \cite{RX12} we proved that if, in addition, each fiber of $f$ in
Theorem A is a topological manifold without boundary and the 
co-dimension is $n$, then $f$ is a fiber bundle projection.
\end{remark}

Theorem A can be viewed as a weak generalization of the local
trivialization of a proper submersion between manifolds.
Recall that a submersion $f$ between two Riemannian manifolds is
said to be an {\it $\epsilon$-Riemannian submersion} if its
differential almost preserves the norm of any horizontal vector,
that is, for any vector $v$ perpendicular to the fibers,
$$e^{-\epsilon}\le\frac{\left|df(v)\right|}{\left|v\right|}\le
e^\epsilon.$$
By definition, an $\epsilon$-Riemannian submersion
is an $e^\epsilon$-LcL. It can be easily checked that 
any smooth $e^\epsilon$-LcL between Riemannian manifolds is an
$e^\epsilon$-Riemannian submersion. 

We shall point it out that an $e^\epsilon$-LcL for $\epsilon>0$ is
weaker than an $\epsilon$-Riemannian submersion or an submetry. 
In the case of an $\epsilon$-Riemannian submersion or an
submetry over Riemannian manifolds, minimal geodesics can be uniquely
lifted to a horizontal curved in the total space and thus give a local
trivialization. For an $e^\epsilon$-LcL, however, the horizontal
curves are not unique any more. Instead, we will apply 
a weaker replacement of horizontal lifting, a neighborhood retraction
$\varphi_p$ to a fiber $f^{-1}(p)$, which is constructed in 
\cite{RX12} (see also \autoref{prop:gradient-retract}, Section 1)
through iterated gradient deformations of distance
functions. If the map is a Riemannian
submersion or a submetry (i.e. $1$-LcL), then
the neighborhood retraction coincides with the local
trivialization map by canonical horizontal liftings. 

Since the proof of Theorem A will be based on local properties, it
still holds if the base space $B$ is a metric space in which each
point admits a neighborhood almost isometric to a ball in $\Bbb R^n$.
In particular, if there is an $(n,\delta)$-strainer at 
a point $p$ in an $n$-dimensional Alexandrov space,
then such bi-Lipschitz homeomorphism of dilation of almost $1$
(depending on $\delta$ and $n$) exists around $p$.
Recall that an {\it$(n,\delta)$-strainer} at $p$ consists of $n$ pairs
of points $(a_i,b_i)_{i=1}^n$ such that the {\it corresponding angles}
(see \cite{BGP92})
$$\tilde\measuredangle a_ipa_j>\frac{\pi}{2}-\delta, \quad
\tilde\measuredangle b_ipb_j>\frac{\pi}{2}-\delta, \quad
\tilde\measuredangle a_ipb_i>\pi-\delta,\quad \tilde\measuredangle
a_ipb_j>\frac{\pi}{2}-\delta.$$ 
And $p$ is said to be $(n,\delta)$-strained. 

\begin{corollary}\label{cor:Hurewicz-fibration}
Given any positive integer $n$, there exists a positive number
$\delta_0(n)$ such that for any $\delta<\delta_0$, there is
$\epsilon>0$ satisfying that if $f$ is a proper
$e^\epsilon$-LcL from a finite-dimensional Alexandrov space $X$ to an
$n$-dimensional Alexandrov space $B$ and each
point of $B$ is $(n,\delta)$-strained,  then $f$ is a Hurewicz
fibration.
\end{corollary}

\begin{remark}
If each fiber of $f$ in \autoref{cor:Hurewicz-fibration} is a
topological manifold without boundary of co-dimension $n$, then by
\autoref{rem:thma}, $f$ is a locally trivial fibration.
\end{remark}

In the case that $f$ is a submetry (i.e., $1$-LcL) it follows from
Perelman's fibration theorem on regular admissible maps (\cite{Pr93})
that $f$ admits a locally trivialization over any point in $B$. This
is because, a submetry $f:X\to B$ satisfies that, for any $p\in B$,
the distance function to the fiber $f^{-1}(p)$ coincides with
$\dist_p\circ f:X\to \Bbb R$ (see \ref{eq:lcl}).
Let $(a_i,b_i)_{i=1}^n$ be an $(n,\delta)$-strainer at $p$.
According to \cite{BGP92}, the map
$\varphi=(\dist_{a_1},\dots,\dist_{a_n}):B_\sigma (p)\to \Bbb R^n$ is
an almost isometry to its image for small $\sigma$, and thus $f$ can
represented by
$$(\dist_{f^{-1}(a_1)},\dots,\dist_{f^{-1}(a_n)})=(\dist_{a_1}\circ
f,\dots, \dist_{a_n}\circ f):X\to \Bbb R^n.$$

In general an $e^\epsilon$-LcL between Alexandrov
spaces fails to satisfy the homotopy lifting property over
singular points. Counterexamples can be constructed by considering
non-free isometric group actions on a Riemannian manifold, where the
quotient space is an Alexandrov space and the projection is a
submetry.

A partial motivation to study an $e^\epsilon$-LcL is that it naturally
arises in some interesting geometry situation. For example, according
to \cite{Ya96} (see also Theorem \ref{thm:Lipschitz-submersion}), if
an Alexandrov space with curvature bounded below is close enough to an
Alexandrov space in Gromov-Hausdorff topology, then an
$e^\epsilon$-LcLs can be constructed over points that are not too
singular. Applications of Theorem A on a convergent
sequence of Alexandrov spaces in Gromov-Hausdorff topology will be
discussed in Section 2, where Yamaguchi's earlier convergence theorem
on Alexandrov spaces is strengthened to a Hurewicz fibration (see
Theorem 2.2) and the same nilpotency results on the fundamental group
of almost nonnegatively curved Alexandrov spaces as Riemannian
manifolds follows from the proof by Kapovitch-Petrunin-Tueschmann in
\cite{KPT10}.

The next main result is about the stability of a
converging sequence of $e^\epsilon$-LcLs. 
Recall that two Hurewicz fibrations $f_i:X_i\to B$ ($i=0,1$) are
fibre-homotopy equivalent if there are fibre-preserving maps $h:X_0
\to X_1$ and $g : X_1 \to X_0$ and fibre-preserving homotopies between
$g\circ h$ and identity $1_{X_0}$, and between $h\circ g$ and
$1_{X_1}$.
We say that Hurewicz fibrations $f_i:X_i\to B_i$ ($i=0,1$) are of the
same homotopy type if there is a homeomorphism $\psi:B_0\to B_1$ such
that $\psi\circ f_0:X_0\to B_1$ is fiber-homotopy equivalent to
$f_1:X_1\to B_1$. For two maps $f_i:A\to B$ $(i=0,1)$ between metric
spaces, we define the distance between $f_0$ and $f_1$ to be
$$d(f_0,f_1)=\sup\{d(f_0(x),f_1(x))\;|\; x\in A\}.$$

\begin{thmb}\label{thm:B}
Let $f_i:X\to B$ $(i=0,1)$ be two $(1.023)$-LcLs from
a finite-dimensional Alexandrov space with curv $\ge \kappa$ to a
Riemannian manifold. Let $r$ be the injectivity radius of $B$. If
$d(f_0,f_1)<\frac{r}{3}$, then $f_0$ and $f_1$ are fiber-homotopy
equivalent. 
\end{thmb}

\begin{remark}
if the two $e^\epsilon$-LcLs in Theorem B are sufficiently close and
the fibers are topological manifolds, then we proved in \cite{RX12}
that they are equivalent as fiber bundles.
\end{remark}

Let $f_i:X_i\to B_i$ ($i=0,1$) be two $(1.023)$-LcLs, where $X_i$
is an $m$-dimensional Alexandrov space with $\curv\ge
\kappa$, $B_i$ is an $n$-dimensional Riemannian manifold
($i=0,1$), and $$d_{GH}((X_0,B_0,f_0),(X_1,B_1,f_1))<\epsilon$$
 in the sense that
 there are $\epsilon$-Gromov-Hausdorff approximations $\varphi:X_0\to
X_1$ and
$\psi:B_0\to B_1$ such that $d(\psi\circ f_0, f_1\circ
\varphi)<\epsilon$.
Recall that an {\it $\epsilon$-Gromov-Hausdorff approximation} is
a (not necessarily continuous) map $\psi:X\to Y$ between metric spaces
such that $| \abs{\psi(x_1)\psi(x_2)}-\abs{x_1x_2}|<\epsilon$
(almost preserving distance) for all $x_1, x_2\in X$ and
$\abs{y\psi(X)} <\epsilon$ (almost dense) for all $y\in Y$. 
According to Perelman's Stability Theorem for Alexandrov spaces with
curvature bounded below (\cite{Pr91}, \cite{Ka07}) and
Cheeger-Gromov's Convergence Theorem (\cite{Ch70}, \cite{GLP81},
\cite{Pr91}, cf. \cite{GW88}, \cite{Pt84}),
there are homeomorphic Gromov-Hausdorff approximations $\phi:X_0\to
X_1$ and $\psi:B_0\to B_1$. By Theorem B, we conclude the following
stability of $e^\epsilon$-LcLs.
\begin{thmc} \label{thm:C}
Let $X_i$ $(i=0,1)$ be an $m$-dimensional Alexandrov space of curv
$\ge \kappa$, $B_i$ $(i=0,1)$ be an $n$-dimensional Riemannian
manifold, and $f_0:X_0\to B_0$ be a $(1.023)$-LcL. Then there
is $\epsilon(X_0,B_0)>0$ such that any $(1.023)$-LcL
$f_1:X_1\to B_1$ satisfying
$d_{GH}((X_0,B_0,f_0),(X_1,B_1,f_1))<\epsilon$
has the same homotopy
type as $f_0$. 
\end{thmc}
\begin{remark}
If, in addition, all fibers of $f_i$ $(i=0,1)$ are closed topological
$(m-n)$-manifolds, then we proved in \cite{RX12} that, there is
$\epsilon_1=\epsilon_1(X_0,B_0,f_0)>0$ such that if
$$d_{GH}((X_0,B_0,f_0),(X_1,B_1,f_1))<\epsilon_1,$$
then $f_1$ is equivalent to $f_0$ as fiber bundles, in the sense that
there are homeomorphic
$\varkappa(\epsilon)$-Gromov-Hausdorff approximations,
$\Psi: X_0\to X_1$, $\Phi:B_0\to B_1$ such that $f_1\circ
\Psi=\Phi\circ f_0$.
Note that a stronger definition of $d_{GH}((X_0,B_0,f_0),$
$(X_1,B_1,f_1))<\epsilon$ is used in \cite{RX12} in the sense
that there are $\epsilon$-Gromov-Hausdorff approximations
$\phi:X_0\to X_1$ and
$\psi:B_0\to B_1$ such that the two maps $\phi\circ f_0$ and
$f_1\circ\psi$ are `fiber-wisely' Hausdorff close:
$$\sup_{\bar x\in
B_1}\{d_H((\phi\circ f_0)^{-1}(\bar x),
(f_1\circ \psi)^{-1}(\bar x))\}<\epsilon,$$
where $d_H$ denotes the Hausdorff distance on subsets in $X_0$.
However, by (\ref{eq:lcl}) these two definitions are equivalent for
LcLs.
\end{remark}

The earlier related result in the smooth category is the fibration
isomorphism/homotopy finiteness of Riemannian submersions under
non-collapsing geometric bounds, which was first proved by Wu
\cite{Wu96}, provided that the fibers are totally geodesic and the
based space is fixed. Later Tapp \cite{Ta00,Ta02} strengthened Wu's
result to general Riemannian submersions. The bundle stability (a
little stronger than the fibration isomorphism finiteness) of
$\epsilon$-Riemannian submersions for small $\epsilon>0$ was proved in
\cite{RX12}. According to \cite{Ka07, Pr91} by Kapovitch and
Perelman, if $f_0$ and $f_1$ in Theorem C are submetries and close
enough to each other, then they are equivalent as fiber
bundles through homeomorphic Gromov-Hausdorff approximations.

Now let us briefly explain the ideas of the proofs of the main
theorems. Recall that a local trivialization for an
$\epsilon$-Riemannian submersion can be directly constructed by the
horizontal liftings of radial minimal geodesics. However, there are no
canonical liftings in case of an $e^\epsilon$-LcL with $\epsilon>0$
due to the lack of regularity.
According to Ferry's result (\cite{Fe78}, see also
\autoref{thm:strong-regular}), the homotopy lifting property
holds for the map in Theorem A if there are controlled homotopy
equivalences between nearby fibers (called {\it strong regular}, see
Section 1.1) and all fibers are abstract neighborhood retracts.
In \cite{RX12} we found a weaker
replacement of horizontal lifting, a neighborhood retraction
$\varphi_p$ to a fiber $f^{-1}(p)$ (see also
\autoref{prop:gradient-retract}, Section 1).
Similar to the horizontal lifting in the smooth case, the neighborhood
retraction $\varphi_p$ of the fiber at $p\in B$  continuously depends
on $p$. Hence the fiber is locally contractible and controlled
homotopy equivalences between nearby fibers can be defined. 
If all fibers are topological manifolds, then
the homotopy equivalences can be approximated by homeomorphisms and
thus $f$ admits a local trivialization (see \cite{RX12}).
Fiber-preserving homotopies for the
composition of maps in Theorem B can be defined similarly.

The remaining of the paper is organized as follows.
In Section 1, we will review some topological results and give 
proofs of Theorem A and Theorem B. In Section 2 we will give  
applications of Theorem A to strengthen convergence theorems for
Alexandrov spaces.

\begin{acknowledgements}
The author would like to thank Xiaochun Rong and Hao Fang for helpful
discussions. The author would also like to thank the University
of Iowa for hospitality and support during a visit in which part of
the work was completed. 
\end{acknowledgements}

\section{Proof of Theorem A and Theorem B}

\subsection{Strong regular maps and semi-concave functions}

Before starting the geometric part of the proofs, we first recall
some topological results. For any Hurewicz fibration $f:X\to Y$, if
$Y$ is path-connected, then by definition the fibers are homotopy
equivalent to each other. In \cite{Fe78} Ferry proved that the
inverse is also true, if the homotopy equivalences between nearby
fibers and the homotopies are under control in the following sense.

A map $f:X\to B$ between metric spaces is said to
be {\it strongly regular} \cite{Fe78} if $f$ is proper and if for each
$p\in B$ and any $\epsilon>0$ there is a $\delta>0$ such that if
$d(p,p')<\delta$,
then there are homotopy equivalences between fibers $\varphi_{pp'}:
f^{-1}(p)\to f^{-1}(p')$, $\varphi_{p'p}: f^{-1}(p')\to f^{-1}(p)$
which togther with the homotopies move points in distance $<\epsilon$.
A topological space $X$ is an {\it absolute neighborhood retract}
(ANR) if there is an embedding of $X$ as a closed subspace of the
Hilbert cube $I^{\infty}$ such that some neighborhood
$N$ of $X$ retracts onto $X$. If $X$ is finite covering dimensional
and locally contractible, then $X$ is an ANR (\cite{Bo55}).

\begin{theorem}[\cite{Fe78}]\label{thm:strong-regular}
If $f:E\to B$ is a strongly regular map onto a complete finite
covering dimensional space $B$ and all fibers are ANRs, then $f$ is a
Hurewicz fibration.
\end{theorem}

According to \autoref{thm:strong-regular}, Theorem A is reduced to
show that an $e^\epsilon$-LcL from an Alexandrov space to a Riemannian
manifold is strongly regular and all its fibers are locally
contractible. 

Now let us recall the gradient flow of a semi-concave function
and the basic property of a LcL that will be frequently used
throughout the paper. Let $X$ be an Alexandrov space with curv
$\ge\kappa$. A
function $f:X\to \Bbb R$ is called {\it semi-concave} (\cite{Petr07})
if for any interior point $p$ in $X$ there is a neighborhood $U$ of
$p$ and a real number $\lambda$ such that for any minimal geodesic
$\gamma(t)$ in $U$, $$f\circ \gamma(t)-\frac{\lambda}{2}t^2$$
is concave. $f$ is called $\lambda$-concave in $U$. The gradient
$\nabla_p f$ as a vector in the tangent cone $T_p$ is well-defined.
For any point $p$ in $X$, there is a unique gradient curve
$\alpha:[0,\infty)\to X$ of $f$
starting at $p$ such that the tangent vector
$\alpha'(t)=\nabla_{\alpha(t)}f$ for any $t\ge 0$. The gradient flow
$\Phi_f^t$ is
well-defined and $e^{\lambda t}$-Lipschitz in $U$ (\cite{Petr07}). Let
$F:\Omega\to X$ be a Lipschitz map from a
metric space, and let $\tau:\Omega\to \Bbb R^+$ be
a Lipschitz function, then the gradient deformation of $F$ with
respect to $f$ is defined to be
$$\Phi_f^{\tau(\omega)}\circ F(\omega):\Omega\to X.$$ 
Now let $f:X\to Y$ be an $e^\epsilon$-LcL between metric spaces. 
For any compact subset $S\subset Y$, let $\operatorname{dist}_S$ be
the distance function to $S$ in $Y$,
$$\operatorname{dist}_S(y)=\abs{y,S}=\inf\{d(y,s): s\in S\}.$$
A basis property of $f$ is that 
the two continuous
functions $\operatorname{dist}_S\circ f$ and
$\operatorname{dist}_{f^{-1}(S)}: X\to
\Bbb R_+$ satisfy (see \cite{RX12})
\addtocounter{theorem}{1}
\begin{equation}\label{eq:lcl}
 e^{-\epsilon}\cdot \operatorname{dist}_S\circ f \le
\operatorname{dist}_{f^{-1}(S)}\le e^\epsilon \cdot
\operatorname{dist}_S\circ f.\tag{1.2}
\end{equation}

\subsection{Proof of Theorem A and Theorem B}
Assume that $X$ is an Alexandrov space with curv $\ge -1$, $B$
is a Riemannian manifold and $f:X\to B$ is an
$e^\epsilon$-LcL. We will try to define controlled homotopy
equivalences between nearby fibers of $f$. As a weaker
replacement of the horizontal lifting of minimal geodesics, a
neighborhood retraction $\varphi_p$ of $f$-fiber over $p\in B$ which
is continuously depending on $p$ was constructed in \cite{RX12}.
Because the proof of Theorem A and Theorem B will rely on the
neighborhood retraction, for reader's convenience, let us recall its
construction in below.

Let $\operatorname{injrad}_B(p)$ denote the injectivity radius of $B$
at $p$. For any point $p\in B$ and $0<r<\operatorname{injrad}_B(p)$,
let $S_r(p)=\partial B_r(p)$ be the metric sphere around $p$ and let
$x$ be any point in $B_r(p)\setminus\{p\}$. It is clear that the
gradient of distance function $\operatorname{dist}_{S_r(p)}$ satisfies
$\left|\nabla_x\dist_{S_r(p)}\right|=1$.  By (\ref{eq:lcl}), an easy
estimate (see Lemma 1.5 in \cite{RX12}) shows that the gradient of
$\operatorname{dist}_{f^{-1}(S_r(p))}$ is bounded by
\addtocounter{theorem}{1}
\begin{equation}\label{eq:gradient-estimate}
 1 -
(e^{2\epsilon}-1) \cdot \frac{2r^2}{\abs{x f^{-1}(p)}\cdot
\abs{x f^{-1} (S_r(p))}} \le \abs{\nabla_x
\operatorname{dist}_{f^{-1}(S_r(p))}}\le 1.\tag{1.3}
\end{equation}
Therefore for sufficient small $\epsilon$ ($e^{\epsilon}\le
1.02368$), points in $f^{-1}(B_{\frac{2r}{3}}(p))$ can be flowed into
$f^{-1}(B_{\frac{r}{3}}(p))$ along gradient curves of
$\operatorname{dist}_{f^{-1}(S_r(p))}$ in a definite time.

\begin{lemma}[Lemma 1.5 in \cite{RX12}]\label{lem:gradient-estimate}
 For any $p\in B$ and $r<\min\{\operatorname{injrad}_p(B),
\frac{1}{2e^\epsilon}\}$, there is a
constant  $C_0(\epsilon)>0$ depending on $\epsilon$ such that for all
$x\in f^{-1}(B_{\frac {2r}3}(p))$,  the gradient
curve $\Phi(t,x)$ of the function
$\operatorname{dist}_{f^{-1}(S_r(p))}$
satisfies
 $$\Phi(x,t)\in f^{-1}(B_{\frac {r}3}(p)),
\qquad t\ge C_0^{-1}\cdot\left(\frac{2}{3}e^{\epsilon} r
-\abs{x f^{-1}(S_r(p))}\right)$$
\end{lemma}

Thus we can define a gradient deformation of
$\operatorname{id}_{f^{-1}(B_{\frac{2r}{3}}(p))}$ which maps
$f^{-1}(B_{\frac{2r}{3}}(p))$ into $f^{-1}(B_{\frac{r}{3}}(p))$ and
fixes $f^{-1}(B_{0.3r}(p))$.
Let
$$T_{p,r}(x)=\max\left\{0,C_0^{-1}\cdot\left(\frac{2}{3}e^{\epsilon}
r -\abs{x f^{-1}(S_r(p))}\right)\right\},$$
and $\Phi_p^{T_{p,r}(x)}(x)=\Phi(x,T_{p,r}(x))$ be the gradient
deformation of $\operatorname{id}_{f^{-1}(B_{\frac{2r}{3}}(p))}$ 
with respect to $\operatorname{dist}_{f^{-1}(S_r(p))}$.
Then by \autoref{lem:gradient-estimate}, for $e^{\epsilon}\le 1.02368$
and $r<\min\{\operatorname{injrad}_p(B),
\frac{1}{2e^\epsilon}\}$, we have
\addtocounter{theorem}{1}
\begin{equation}\label{eq:gradient-flow}
 \begin{cases}
 \Phi_p^{T_{p,r}(x)}(x)\in f^{-1}(B_{\frac {r}3}(p)), &\forall\; x\in
f^{-1}(B_{\frac{2r}{3}}(p)),\\
 T_{p,r}(x)=0, &\forall\; x\in
f^{-1}(B_{0.3 r}(p)).\tag{1.5}
   \end{cases}
\end{equation}
Because the gradient curves are stable as function converges
(\cite{Petr07}), it is proved in \cite{RX12} that
$\Phi_p^{T_{p,r}(x)}(x)$ is continuous both in $p$ and $x$, that is,
$$ \Psi:\bigcup_{p\in B} \{p\}\times
f^{-1}(B_{\frac{2r}{3}}(p))\subset B\times X
\to X,\quad \Psi(p,x)=\Phi_p^{T_r(x)}(x)$$
is a continuous map.

Repeating the construction above for the sequence
$\{r_i=\frac{r}{2^i}\}_{i=0,1,2,\cdots}$
and let $\Phi_{p,i}^{T_{p,i}}(x)=\Phi_{p,i}(x,T_{p,r_i}(x))$ be the
gradient curves of $\operatorname{dist}_{f^{-1}(S_{r_i}(p))}$ at $x$
with time $T_{p,r_i}(x)$. By (\ref{eq:gradient-flow}),
$\Phi_{p,i}^{T_{p,i}}:f^{-1}(B_{r_i}(p))\to X$
takes $f^{-1}(B_{\frac{2}{3}\cdot \frac{r}{2^i}}(p))$ into 
$f^{-1}(B_{\frac{1}{3}\cdot\frac{r}{2^{i+1}}}(p)),$
and $$\left .\Phi_{p,i}^{T_{p,i}}
\right |_{f^{-1}(B_{0.3\frac{r}{2^i}}(p))} = \operatorname{id}.$$
Hence the iterated gradient deformations
$$\Phi_{p,i}^{T_{p,i}}\circ \Phi_{p,i-1}^{T_{p,i-1}}\circ \cdots
\Phi_{p,0}^{T_{p,0}}$$
is well-defined on $f^{-1}(B_{\frac{2r}{3}}(p))$ and its restriction
on $f^{-1}(B_{0.3\frac{r}{2^i}}(p))$ is identity.
Because $$T_{p,r_i}(x)\le \frac{r}{2^{i-1}}\cdot
\frac{e^\epsilon}{3}\cdot C_0^{-1},$$ it can be
directly verified that the sequence of maps
\begin{gather*}
 \Psi_i:\bigcup_{p\in B} \{p\}\times f^{-1}(B_{\frac{2r}{3}}(p))
\to X,\\
\Psi_i(p,x)=\Phi_{p,i}^{T_{p,i}}\circ \Phi_{p,i-1}^{T_{p,i-1}}\circ
\cdots \Phi_{p,0}^{T_{p,0}}(x)
\end{gather*}
uniformly converges. The limit $\varphi_p(x)=\displaystyle\lim_{i\to
\infty}\Psi_i(p,x)$ gives a retraction
from the neighborhood $f^{-1}(B_{\frac{2r}{3}}(p))$ to $f^{-1}(p)$,
which is continuous both in $p$ and $x$. 

\begin{proposition}[Proposition 1.6 in \cite{RX12}]
\label{prop:gradient-retract}
For any $0<r<\injrad(B)$, there is a map
$\varphi_p(x)$ from a neighborhood $f^{-1}(B_{\frac{2r}{3}}(p))$ to
the fiber $f^{-1}(p)$ such that  
$$\varphi_p(x): \bigcup_{p\in B} \{p\}\times
f^{-1}(B_{\frac{2r}{3}}(p))\to X,$$
is continuous both in $p$ and $x$, and 
satisfies 
\begin{enumerate}
 \item[(\ref{prop:gradient-retract}.1)] $\varphi_p(x)=x$ for any $x\in
f^{-1}(p)$, and 
 \item[(\ref{prop:gradient-retract}.2)] $\abs{x \varphi_p(x)} \le
2C_1r,
$ for some constant
$C_1(\epsilon)$ depending only on $\epsilon$.
\end{enumerate}
\end{proposition}

\begin{remark}
Since the Lipschitz constant of each $\Phi_i$ goes to infinity with
the same order as $r_i^{-1}$, there is no Lipschitz
control on the limit map $\varphi_p$. Therefore no local control on
the intrinsic distance in the fiber follows from
\autoref{prop:gradient-retract}.
\end{remark}

We are now ready to prove Theorem A and Theorem B.
\begin{proof}[Proof of Theorem A]
Up to a rescaling we assume that the lower curvature bound of $X$ is
$-1$.
By \autoref{thm:strong-regular}, it suffices to
show that $f$ is strong regular and any fiber is an ANRs. 
For any $p,q\in B$ with small distance
$0<\abs{pq}<\frac{1}{2}\min\{\operatorname{injrad}_p(B),
\frac{1}{2e^\epsilon}\}$, let $\rho=2\abs{pq}$. By
\autoref{prop:gradient-retract}, there are neighborhood
retractions $\varphi_p:
f^{-1}(B_{\frac{2\rho}{3}}(p))\to f^{-1}(p)$
and $\varphi_q: f^{-1}(B_{\frac{2\rho}{3}}(q))\to f^{-1}(q)$
around $f^{-1}(p)$ and $f^{-1}(q)$ respectively.
Then the homotopy equivalences between fibers can be chosen to be
$\left.\varphi_p\right|_{f^{-1}(q)} :f^{-1}(q)\to f^{-1}(p)$
and $\left.\varphi_q\right|_{f^{-1}(p)}:f^{-1}(p)\to f^{-1}(q)$, and
the homotopies are
$H_t=\varphi_p\circ \varphi_{\gamma(t)}:f^{-1}(p)\to f^{-1}(p)$ and
$K_t=\varphi_q\circ \varphi_{\gamma(1-t)}:f^{-1}(q)\to f^{-1}(q)$,
where $\gamma:[0,1]\to B$ is a minimal geodesic from $p$ to $q$.
By (\ref{prop:gradient-retract}.2), $\abs{H_t(x)x}\le 4C_1\rho$ and
$\abs{K_t(x)x}\le 4C_1\rho$. Therefore $f$ is strongly regular.

From \cite{Pr91}, an Alexandrov space with curvature bounded below is
locally contractible. For $x\in
f^{-1}(p)$, let $U_x\ni x$ be a contractible neighborhood around
$x$ and $H_t:U_x\to U_x$ be the homotopy from
$\operatorname{id}_{U_x}$ to the retraction $r:U_x\to \{x\}$ such
that $H_t(x)=x$. Then $\varphi_p\circ H_t$ is a homotopy
from $\operatorname{id}_{U_x\cap f^{-1}(p)}$ to the retraction $r
:U_x\cap f^{-1}(p)\to \{x\}$. Therefore $f^{-1}(p)$ is locally
contractible and thus an absolute neighborhood retract.
\end{proof}

\begin{proof}[Proof of Theorem B]
Let $f_0,f_1:X\to B$ be the $(1.023)$-LcLs from Alexandrov space $X$
to
Riemannian manifold $B$ in Theorem D. We now construct 
fiber-preserving maps $h,g:X\to X$ and fiber-preserving homotopies 
$g\circ h$ to the identity $1_A$ and from $h\circ g$ to $1_A$ as
follows.

For any point $x\in X$, let $p=f_0(x)\in B$, let $F_0(p)$ be
the fiber $f^{-1}_0(p)$ and $F_1(p)=f^{-1}_1(p)$. Because
$d(f_0,f_1)<\frac{r}{3}$, by (\ref{eq:lcl}), $F_0(p)$ lies in the
$1.023\frac{r}{3}$-tubular neighborhood of $F_1(p)$. Let $\varphi_p$
be the neighborhood retraction of $F_1(p)$ in
\autoref{prop:gradient-retract} with respect to $f_1$,  we define
$h:X\to X$ by $h(x)=\varphi_p(x)=\varphi_{f_0(x)}(x)$. Then the
continuous map $h:X\to X$ is globally defined and maps all fibers of
$f_0$ into that of $f_1$. Similarly we define $g:X\to X$ through the
neighborhood retraction of $f_0$-fibers such that $f_0\circ
g=f_1$, where $g(x)=\psi_{f_1(x)}(x)$ and $\psi_{q}$ is the
neighborhood retraction of $f_0^{-1}(f_1(x))$ with respect to $f_0$.

Let $\gamma:[0,1]\to B$ be the minimal
geodesic from $\gamma(0)=f_0(x)$ to $\gamma(1)=f_1(x)$, we define the
fiber-preserving homotopy $H_t:X\to X$ by
$H_t(x)=\psi_{f_0(x)}\circ \varphi_{\gamma(t)}(x)$. Then
$H:[0,1]\times X\to X$ is a continuous map such that $H_0=g\circ h$
and $H_1=1_{X}$. A fiber-preserving homotopy from $h\circ g$
to the identity $1_A$ can be defined similarly.
\end{proof}

\section{Convergence theorems of Alexandrov spaces}
In this section we will apply Theorem A to
strengthen Yamaguchi's Lipschitz submersion theorem on Alexandrov
spaces. A convergent sequence of Riemannian manifolds or
Alexandrov spaces $X_i$ in Gromov-Hausdorff topology is called {\it
collapsing} if the dimension of limit
space is strictly less than that of $X_i$. 
 for $\delta>0$, the {\it $\delta$-strain radius} \cite{Ya96} at a
point $p$ in an $n$-dimensional Alexandrov space $X$ is
defined to be $$\sup \{r\; | \text{ there exists an
$(n,\delta)$-strainer
at $p$ of length $r$} \}.$$ 
A Lipschitz submersion
theorem was proved by Yamaguchi \cite{Ya96} in the case that all
points in the limit space is $(n,\delta)$-strained.
\begin{theorem}[Lipschitz submersion Theorem
\cite{Ya96}]\label{thm:Lipschitz-submersion}
For any dimension $n$ and positive number $\mu_0$, there exist
positive numbers $\delta(n)$ and $\epsilon(n,\mu_0)$ such that
for any $m$-dimensional Alexandrov space $Y$ with curv $\ge -1$ and
any
$n$-dimensional Alexandrov space $X$ with curv $\ge -1$, if
\begin{enumerate}
 \item the $\delta$-$\operatorname{strain\; radius}$ at any point in
$X$ $\ge \mu_0$, and
 \item  the Gromov-Hausdorff distance between $X$ and $Y$ $\le
\epsilon$, then
\end{enumerate}
there exists a $(1+\tau(\delta,\epsilon))$-LcL $f : Y \to
X$ which is a $\tau(\delta,\epsilon)$-almost Lipschitz 
submersion.

Here, $\varkappa(\delta,\epsilon)$ denotes a positive
constant depending on $n$, $\delta$ and $\epsilon$ satisfying
$\tau(\delta,\epsilon)\to 0$ as $\delta,\epsilon\to 0$.
\end{theorem}

Yamaguchi conjectured in \cite{Ya96} that the map $f:Y\to X$ in
\autoref{thm:Lipschitz-submersion} is a locally trivial fibration.
By Theorem A, we are able to conclude that it
is a Hurewicz fibration, and the
conjecture is true in the case that the fibers are topological
manifolds. 

\begin{theorem}
For sufficient small $\delta$, the almost Lipschitz 
submersion $f:Y\to X$ in
\autoref{thm:Lipschitz-submersion} is a Hurewicz fibration; and $f$
is a locally trivial fibration in the case that
every $f$-fiber is a topological manifold
without boundary of co-dimension $n$.
\end{theorem}

There are similar fibration theorems for Riemannian
manifolds which are fundamental in many applications.
It was first proved by
Fukaya in \cite{Fu87} and strengthened later by Cheeger-Gromov-Fukaya
in \cite{CFG92} that if a sequence of Riemannian manifold $M_i$ is
collapsing under a two-sided sectional curvature bound,
then the Gromov-Hausdorff approximation can be approximated by a
global singular fibration from $M_i$ to the limit space, whose fibers
are almost flat manifolds. 
Yamaguchi proved in \cite{Ya91} that if a Riemannian
manifold $M$ with sectional curvature bounded below is sufficient
close to a Riemannian manifold $N$ in Gromov-Hausdorff distance, then
there is a $C^1$ $\epsilon$-Riemannian submersion $f:M\to N$,
whose fiber is connected and has almost nonnegative curvature in a
generalized sense defined in Lemma 5.3 of \cite{Ya91} (cf. Definition
1.0.1 in \cite{KPT10}).  By Kapovitch-Petrunin-Tueschmann's
work \cite{KPT10} the fundamental group of the fiber admits a
nilpotent subgroup of uniformly bounded index.

Due to Theorem 2.2 one is able to talk about the homotopy fiber.
By similar arguments as the smooth case, the same nilpotency on the
fundamental group of the homotopy fiber of $f$ in
\autoref{thm:Lipschitz-submersion} follows from
Theorem A and the proof in \cite{KPT10} by
Kapovitch-Petrunin-Tueschmann.

\begin{corollary}
There is a universal constant depending on $(m-n)$ such that
the fundamental group of the homotopy $F$ in
\autoref{thm:Lipschitz-submersion} contains a nilpotent
subgroup whose index $\leq C$.
\end{corollary}

Almost nonnegatively curved Alexandrov spaces are those space whose
diameter is uniformly bounded and the lower curvature bound $\kappa$
is almost nonnegative. Or equivalently, the scaling invariant
$\operatorname{diam}^2\cdot \kappa$ is almost nonnegative.
After scaling to a uniform lower curvature bound $-1$, they are
close to a point in Gromov-Hausdorff topology. As a special case, the
almost nilpotency for almost nonnegatively curved Alexandrov spaces
follows.

\begin{theorem}
For any positive integer $n$, there are constants $\epsilon(n)>0$ and
$C(n)>0$ such that for any closed $n$-dimensional Alexandrov space
$M$ with curv $\geq \kappa$, if
$\operatorname{diam}(M)^2\cdot \kappa>-\epsilon(n)$, then $\pi_1(M)$
contains a nilpotent subgroup of index $\le C(n)$.
\end{theorem}

\begin{remark}
The earlier version of Theorem 2.4 for Alexandrov spaces
was proved by Yamaguchi in \cite{Ya96} without a uniform bound on the
index of the nilpotent subgroup, where the proof was based the
Lipschitz submersion theorem (\autoref{thm:Lipschitz-submersion})
and arguments in \cite{FY92}.
\end{remark}

Theorem 2.4 extends the same property of almost nonnegatively curved
manifolds. It was first proved by Fukaya and Yamaguchi \cite{FY92}
that the fundamental group of any manifold of almost nonnegative
sectional curvature contains a nilpotent subgroup of finite index,
and later Kapovitch, Petrunin and Tueschmann \cite{KPT10} proved
that the index of a nilpotent subgroup is uniformly bounded
(depending only on the manifold's dimension). 
The same conclusion for manifolds of almost nonnegative Ricci
curvature is a conjecture of Gromov, and it was proved recently by
Kapovitch and Wilking in \cite{KW11} (cf. \cite{CC96}).

Recall that Yamaguchi's fibration theorem \cite{Ya91} for Riemannian
manifolds provided a fundamental tool in constructing a finite
descending normal tower of the fundamental group in earlier proofs
both in \cite{FY92} and \cite{KPT10}.
By the fibration in Theorem 2.2, we still are able to conclude
the same structure for a collapsing sequence of Alexandrov spaces.
Moreover, the estimates in \cite{KPT10} were proved for general
Alexandrov spaces.
Because the proof of Theorem 2.4 and Corollary 2.3
would be the same as that in \cite{KPT10}, we omit it here.

Another direct corollary of Theorem 2.2 is a long exact sequence
arising from
the fibration.
\begin{corollary}\label{cor:collapsing-long-exact-sequence}
 Assume a sequence of Alexandrov spaces $X_i$ with curv $\ge\kappa$
collapsing to an $n$-dimensional Alexandrov $X$ and the limit space
$X$ has only weak singularities. Then for all large $i$, the almost
Lipschitz submersion $f_i:X_i\to X$ induces an isomorphism
$f_*:\pi_l(X_i,F_i,x_i)\to \pi_l(X,x)$ for all $l\ge
1$, where $x\in X$ and $x_i\in F_i=f_i^{-1}(x)$.
 Hence there is a long exact sequence
 $$\cdots\to\pi_l(F_i,x_i)\to
\pi_l(X_i,x_i)\overset{f_*}{\to}\pi_l(X,x)\to\pi_{l-1}(F_i,
x_i)\to\cdots\to\pi_1(X,x)\to 0.$$
\end{corollary}

\begin{remark}
In \cite{Pr97} Perelman concluded the same long exact sequence 
under a weaker condition that the limit $X$ has no proper extremal
subsets,
where $F_i$ is a regular fiber (i.e., the fiber of a lifting map to
$X_i$ of regular admissible maps locally defined on the limit $X$ to
$\Bbb R^n$
\cite{Pr97}). 
It is easy to see that the regular fiber is homotopy equivalent to
$f_i^{-1}(x)$ in Corollary \ref{cor:collapsing-long-exact-sequence}.
\end{remark}

\bibliographystyle{plain}
\bibliography{almost_submetry_references}

\begin{thebibliography}{10}

\bibitem{Bo55}
K.~Borsuk.
\newblock On some metrizations of the hyperspace of compact sets.
\newblock {\em Fund. Math.}, 41:168--201, 1955.

\bibitem{BGP92}
Y.~Burago, M.~Gromov, and G.~Perelman.
\newblock A.d. alexandrov spaces with curvature bounded below.
\newblock {\em Uspekhi Mat. Nauk}, 47(2(284)):3--51, 1992.

\bibitem{Ch70}
J.~Cheeger.
\newblock Finiteness theorems for riemannian manifolds.
\newblock {\em Amer. J. Math.}, 92:61--75, 1970.

\bibitem{CC96}
J.~Cheeger and T.~H. Colding.
\newblock Lower bounds on ricci curvature and the almost rigidity of warped
  products.
\newblock {\em Ann. of Math.}, 144(1):189--237, 1996.

\bibitem{CFG92}
J.~Cheeger, K.~Fukaya, and M.~Gromov.
\newblock Nilpotent structures and invariant metrics on collapsed manifolds.
\newblock {\em J. A. M. S.}, 5:327--372, 1992.

\bibitem{Fe78}
S.~Ferry.
\newblock Strongly regular mappings with compact anr fibers are hurewicz
  fiberings.
\newblock {\em Pacific J. Math.}, 75(2):373--382, 1978.

\bibitem{Fu87}
K.~Fukaya.
\newblock Collapsing of riemannian manifolds to ones of lower dimensions.
\newblock {\em J. Diff. Geom.}, 25:139--156, 1987.

\bibitem{FY92}
K.~Fukaya and T.~Yamaguchi.
\newblock The fundamental groups of almost non-negatively curved manifolds.
\newblock {\em Ann. of Math.}, 136:253--333, 1992.

\bibitem{GW88}
R.~E. Green and H.~Wu.
\newblock Lipschitz convergence of riemannian manifolds.
\newblock {\em Pacific J. Math.}, 131:119--141, 1988.

\bibitem{GLP81}
M.~Gromov, J.~Lafontaine, and P.~Pansu.
\newblock {\em Structures M\'{e}triques pour les Vari\'{e}t\'{e}s
  Riemannienes}.
\newblock CedicFernand, Paris, 1981.

\bibitem{Ka07}
V.~Kapovitch.
\newblock Perelman's stability theorem.
\newblock In {\em Surveys of Differential Geometry, Metric and Comparison
  Geometry}, XI, pages 103--136. Int. Press, Somerville, 2007.

\bibitem{KPT10}
V.~Kapovitch, A.~Petrunin, and W.~Tuschmann.
\newblock Nilpotency, almost nonnegative curvature, and the gradient flow on
  alexandrov spaces.
\newblock {\em Ann. of Math.}, 171:343--373, 2010.

\bibitem{KW11}
V.~Kapovitch and B.~Wilking.
\newblock Structure of fundamental groups of manifolds with ricci curvature
  bounded below.
\newblock {\em Preprint, arXiv:1105.5955}, 2011.

\bibitem{Pr91}
G.~Perelman.
\newblock Alexandrov spaces with curvatures bounded from below ii.
\newblock {\em Preprint}, 1991.

\bibitem{Pr93}
G.~Perelman.
\newblock Elements of morse theory on aleksandrov spaces.
\newblock {\em St. Petersburg Math. J.}, 5:205--213, 1993.

\bibitem{Pr97}
G.~Perelman.
\newblock Collapsing with no proper extremal subsets.
\newblock In K.~Grove and P.~Petersen, editors, {\em Comparison Geometry},
  volume~30 of {\em MSRI Books}, pages 149--155. Cambridge Univ. Press,
  Cambridge, 1997.

\bibitem{Pt84}
S.~Peters.
\newblock Cheeger's finiteness theorem for diffeomorphism classes of riemannian
  manifolds.
\newblock {\em J. Reine Angew. Math.}, 349:77--82, 1984.

\bibitem{Petr07}
A.~Petrunin.
\newblock Semiconcave functions in alexandrov's geometry.
\newblock In J.~Cheeger and K.~Grove, editors, {\em Metric and Comparison
  Geometry}, volume~XI of {\em Surveys in Differential Geometry}, pages
  137--202. Int. Press, Somerville, 2007.

\bibitem{RX12}
X.~Rong and S.~Xu.
\newblock Stability of $e^\epsilon$-lipschitz and co-lipschitz maps in
  gromov-hausdorff topology.
\newblock {\em Advances in Mathematics}, 231:774--797, 1 October 2012.

\bibitem{Ta00}
K.~Tapp.
\newblock Bounded riemannian submersions.
\newblock {\em Indiana Univ. Math. J.}, 49(2):637--654, 2000.

\bibitem{Ta02}
K.~Tapp.
\newblock Finiteness theorems for submersions and souls.
\newblock {\em Proc. Amer. Math. Soc.}, 130(6):1809--1817, 2002.

\bibitem{Wu96}
J.~Y. Wu.
\newblock A parametrized geometric finiteness theorem.
\newblock {\em Indiana Univ. Math. J.}, 45(2):511--528, 1996.

\bibitem{Ya91}
T.~Yamaguchi.
\newblock Collapsing and pinching under a lower curvature bound.
\newblock {\em Ann. of Math.}, 133:317--357, 1991.

\bibitem{Ya96}
T.~Yamaguchi.
\newblock A convergence theorem in the geometry of alexandrov spaces, 1996.

\end{thebibliography}

\end{document}